\documentclass[11pt]{article} 
\usepackage{amsmath}
\usepackage{makeidx}
\usepackage{graphicx}
\usepackage{epsfig}
\usepackage{amsfonts}

\newtheorem{theorem}{Theorem}

\newtheorem{lemma}{Lemma}
\newtheorem{proposition}[lemma]{Proposition}

\def\intl#1{\hbox{$\backslash$}\kern-10pt\int_{#1}}
\def\weakto{\rightharpoonup}
\def\div{\mbox {div }}

\begin{document}
\title{\vspace{-3cm} 
Partial
 Regularity of Navier-Stokes Equations}
\author{Lihe Wang}
%
\date{}
\maketitle
\begin{abstract}
We prove, with a more geometric approach,  that the solutions  to the  Navier-Stokes equations are regular up to a set of Hausdorff dimension 1.
The main tool for the proof is a new compactness lemma and the monotonicity property of harmonic functions.
\end{abstract}
\section {Introduction }

In the seminal work by L. Caffarelli, R. Kohn and Nirenberg \cite{CKN}, they proved that the solution of Navier-Stokes equations
\begin{equation}\label{NS}
\left\{\begin{array}{rl}
u_t+u\cdot \nabla u-\Delta u+\nabla P&=0\\
\div u&=0
\end{array}\right.
\end{equation} is smooth away from a set of zero one Hausdorff dimensional measure in the parabolic metric.
The original proof in \cite{CKN} consists of a series of  delicate iterations which have been simplified to some extend,  by Fanghua Lin \cite{Lin},O. A. Ladyzhenskaya and G. A. Seregin \cite{LS}, and many others.
These monumental  iterations are nevertheless very much specialized to Navier-Stokes equations.
This work is an attempt to provide a  geometric framework and to  show how the regularity is propagated cross different scales.  We will use the  classical monotonicity inequality for harmonic functions as an essential tool for  controlling small scales from the unit scale.  This makes our computation geometrical and apparent. 
We also introduce a compactness lemma which is suitable to handle the  lack of compactness in the pressure in the time variable, 
This shall have wide applications in other situations.

\section{The Monotonicity formula}
We use the standard notations for parabolic equations as $B_r=\{x:|x|<r\}$ and $Q_r= B_r\times (-r^2, 0]$. We also use $\displaystyle\intl{A} f$ for the average of $f$ in $A$ as
$\frac 1{|A|}\int_{A}f$  where the integration is Lebesgue integral and $|A|$ is its measure in its dimension.

First we recall the classical monotonicity inequality for harmonic functions.

\begin{lemma} \label{harmonic}Suppose $u$ is a harmonic function in $B_1$, then  for any $1\leq p<\infty$, the following functions are monotone increasing in $r$: 
 \begin{equation}\label{mono}\intl{B_r}|u|^p dx,\\
 \end{equation}
 \begin{equation}\label{boundarymono}
 \intl{\partial B_r}|u|^p d\sigma.\\
 \end{equation}

\end{lemma}

We prove several perturbations of the above lemma.

\begin{lemma}[Inhomogeneous Monotonicity]\label{inmono}
 
 (a) For any $1\leq q<n$ and $1\leq p\leq p^*:=\frac {nq}{n-q}$, there are universal constants $C_1$ and $C_2$ such that for any    weak solution $u$  in $B_1\subset \mathbb{R}^n$
 \begin{eqnarray*}
\Delta u&=& \div F, 
\end{eqnarray*}
we have  for any  $ 0<r\leq 1 $: 
 \begin{equation}\intl{B_r}|u|^p dx\leq C_1\intl{B_1}|u|^p dx+\frac {C_2}{r^n}\left(\intl{B_1} |F|^q\right)^{\frac {p}q}.\end{equation}
 
 (b)For any $1< p<\infty$, there are universal constants $C_1$ and $C_2$ such that for any  distribution solution $u$  in $B_1\subset R^n$ of
 \begin{eqnarray*}
\Delta u&=& \div^2 G:=\sum_{i,j=1}^n(G_{ij})_{x_ix_j}, 
\end{eqnarray*}
we have  for any   $ 0<r\leq 1 $: 
 \begin{equation}
\label{CZ}
\intl{B_r}|u|^p dx\leq C_1\intl{B_1}|u|^p dx+\frac {C_2}{r^n}\intl{B_1} |G|^p.\end{equation}
\end{lemma}
{\bf Proof.} The proofs for (a)  and (b) are completely parallel so we present a proof of (b) only. Also (a) is a consequence of (b) by  a proper choice of $G$ using  the Sobolev embedding. Now we prove (\ref{CZ}).  Since the  inequality  (\ref{CZ}) is clearly true for $\frac 12\leq r\leq 1$ if $C_1\geq 2^n$, we prove it  only for $0<r<\frac 12.$
First, by Fubini's theorem, we can find $r_0\in [\frac 12,\frac 34]$ so that 
 \begin{eqnarray*}
\frac{1}{r_0^{n-1}}\int_{\partial B_{r_0}}|u|^p d\sigma &
\leq &C\intl{B_1}|u|^p dx.
\end{eqnarray*}
Now let $h$ be the harmonic function
 $$\left\{\begin{array}{rcl}
\Delta h&=&0,  \mbox{ in } B_{r_0}\\[2mm]
h&=& u, \mbox{ on  }\partial B_{r_0}.
\end{array}\right.
$$
Then we have, with $r\leq \frac 12\leq r_0$,  via Lemma \ref{harmonic},
\begin{eqnarray}
\label{1l}\intl{B_{r}}|h|^p dx &\leq &\displaystyle\intl{B_{r_0}}|h|^p dx\notag\\
&=&\frac{n}{r_0}\displaystyle\int_0^{r_0}\left(\intl{\partial B_\rho}|h|^p\rho^{n-1} d\sigma \right)d\rho\notag\\
&\leq &\frac{1}{r_0}\intl {\partial B_{r_0}}|h|^p dx\notag\\
&\leq &C\displaystyle\intl{B_1}|u|^p dx\label{boundary}
\end{eqnarray}
 and also that 
 $$\left\{\begin{array}{rcl}
\Delta (u-h)&=&div^2 G,  \mbox{ in } B_{r_0}\\[2mm]
u-h&=& 0, \mbox{ on  }\partial B_{r_0},
\end{array}\right.
$$
which yields, with Caldron-Zygmund estimates, that
 \begin{eqnarray*}\label{2l}\intl{B_{r_0}}|u-h|^{p} dx&\leq &C \intl{B_{r_0}}|G|^{p} dx\\
&\leq&\frac{C}{r_0^n}\intl{B_{1}}|G|^{p} dx\\
&\leq &2^nC\intl{B_{1}}|G|^{p} dx.
\end{eqnarray*}
Therefore,
 \begin{eqnarray*}
\intl{B_r}|u|^p dx&\leq& 2^{p-1}\left(\intl{B_r}|h|^p dx+\intl{B_r}|h-u|^p dx\right)\\
&\leq & 2^{p-1}\left(\intl{B_{r_0}}|h|^p dx+\frac {r_0^n}{r^n}\intl{B_{r_0}}|h-u|^p dx\right)\\
&\leq & C\intl{B_{1}}|u|^p dx+\frac {C}{r^n}\intl{B_1}|G|^pdx.
\end{eqnarray*}

The following  is a variant of the above monotonicity formula which works only in 3-dimension.
\begin{lemma}[Interpolated Monotonicity]\label{linter}
 There are universal constants $C_1$ and $C_2$ such that for any    $u$  in $H^1(B_1)$ for $ B_1\subset  R^3$,  we have  for   $ 0<r\leq 1 $: 
 \begin{equation}\label{ininter}
 \intl{B_r}|u|^3 dx\leq C_1\intl{B_1}|u|^3 dx+\frac {C_2}{r^3}\left(\int_{B_\frac 34}|u|^2\right)^\frac 12\int_{B_1} |Du|^2.\end{equation}
\end{lemma}
{\bf Proof.} While the  inequality (\ref{ininter})  is clearly true for $\frac 12\leq r\leq 1$ by taking $C_1\geq 2^n$, we consider $0<r<\frac 12.$
We first write
 \begin{eqnarray*}
\Delta u&=& \div F, 
\end{eqnarray*}
with $F=\nabla u$. 
Then by Fubini theorem and the mean value inequality, we find $r_0\in [\frac 12,\frac 58]$ so that  
 \begin{eqnarray*}
\frac{1}{r_0^{2}}\int_{\partial B_{r_0}}|u|^3 d\sigma &\leq &C\intl{B_1}|u|^3 dx,\\
\frac{1}{r_0^{2}}\int_{\partial B_{r_0}}|u|^2 d\sigma &\leq &C\intl{B_\frac34}|u|^2 dx.
\end{eqnarray*}
As before we  let $h$ be the harmonic function
$$\left\{ \begin{array}{rcl}
\Delta h&=&0,  \mbox{ in } B_{r_0}\\
h&=& u, \mbox{ on  }\partial B_{r_0}.
\end{array}\right.
$$
With a simple computation as (\ref{boundary}), we have 
\begin{equation}
\intl{B_{r_0}}|h|^2\leq C\intl{\partial B_{r_0}} |u|^2\leq C\intl{B_\frac34}|u|^2.
\label{h2}
\end{equation}
Furthermore from the monotonicity formula in Lemma \ref{harmonic}  for any $r\leq \frac 12\leq r_0$,
\begin{equation}\label{1}\intl{B_{r}}|h|^3 dx\leq \intl{B_{r_0}}|h|^3 dx\leq C\intl{\partial B_{r_0}}|h|^3 dx\leq C\intl{B_1}|u|^3 dx.\end{equation}
We also have
 $$
 \left\{\begin{array}{rcl}
\Delta (u-h)&=&\div F,  \mbox{ in } B_{r_0}\\
u-h&=& 0, \mbox{ on  }\partial B_{r_0},
\end{array}\right.
$$
which yields that, with the Sobolev embedding and energy estimates,
 \begin{eqnarray*}\label{2}\intl{B_{r_0}}|u-h|^{6} dx&\leq  & C {r_0^2}\left(\intl{B_{r_0}}|D(u-h)|^{2} dx\right)^{3}\\
&\leq &C \frac{1}{r_0}\left(\int_{B_{r_0}}|F|^{2} dx\right)^{3}\\
&\leq &2^7C\left(\int_{B_{1}}|F|^{2} dx\right)^{3}.
\end{eqnarray*}
Similarly with Poincar\'{e} inequality,
$$
\intl{B_{r_0}}|u-h|^2\leq Cr^2_0\intl{B_{r_0}}|D(u-h)|^2\leq Cr^2_0\intl{B_{r_0}}|F|^2\leq 4C\intl{B_1}|Du|^2.
$$
Finally we derive (\ref{ininter}) as
 \begin{eqnarray*}
\intl{B_r}|u|^3 dx&\leq& 4\left(\intl{B_r}|h|^3 dx+\intl{B_r}|h-u|^3 dx\right)\\
&\leq & 4\left(\intl{B_{r_0}}|h|^3 dx+\frac {r_0^3}{r^3}\intl{B_{r_0}}|h-u|^3 dx\right)\\
&\leq & C\intl{B_{1}}|u|^3 dx+\frac {C}{r^3}\left(\int_{B_{r_0}}|u-h|^2\right)^\frac 34\left(\int_{B_{r_0}}|h-u|^6\right)^{\frac 14}\\
&\leq &C\intl{B_{1}}|u|^3 dx+\frac {C}{r^3}\left(\int_{B_{r_0}}|u-h|^2\right)^\frac 34\left(\int_{B_{1}}|Du|^2\right)^{\frac 34}\\
&\leq &C\intl{B_{1}}|u|^3 dx+\frac {C}{r^3}\left(\int_{B_{r_0}}|u-h|^2\right)^\frac 12\left(\int_{B_{1}}|Du|^2\right),\end{eqnarray*}
where we have used the previous inequality. Inequality (\ref{ininter}) follows from the above inequality and the fact that,
$$
\left(\int_{B_{r_0}}|u-h|^2\right)^\frac 12\leq \left(\int_{B_{r_0}}|u|^2\right)^\frac 12+\left(\int_{B_{r_0}}|h|^2\right)^\frac 12\leq C\left(\int_{B_\frac34} |u|^2\right)^\frac 12,
$$
where we have used  (\ref{h2}).

\section{Suitable weak solutions}
From now on, we will consider (1) in 3-dimension. 
We recall that if $\{u, P\}$ is a suitable weak solution  as in  (\ref{NS}) if    $u$ is in the space $ L^2((0, T), H^1)\cap H^1((0, T), H^{-1})$ and that $P\in L^\frac 32$  so that 

(a) $\{u, P\} $ solves (\ref{NS}) in the weak sense  while the term $u\cdot \nabla u$ is understood as $\partial_{i}(u^i u^j) $ weakly and that

(b) for any $\varphi\geq 0$ in $ C^\infty_0(Q_1)$, then
\begin{equation}
\begin{array}{rl}\displaystyle\int_{Q_1\cap\{t\}}|u|^2\varphi dx+2\int_0^t\int_{B_1}|Du|^2\varphi dxdt&\\
&\hspace{-2in}\displaystyle\leq \int_0^t\int_{B_1}\left[|u|^2(\varphi_t+\Delta \varphi)+(|u|^2+2P)u\cdot D\varphi \right]dxdt.
\end{array}\label{energy}
\end{equation}

We prove the following well known facts for the reader's convenience.

\begin{lemma}[Local energy estimate]\label{local-energy} 
Suppose  $\{u, P\}$ is a suitable weak solution of (\ref{NS}) in $Q_1$. Then 
$$\displaystyle\sup_{t\in (-\frac {9}{16}, 0]}\int_{B_\frac 34 }|u(\cdot, t)|^2+\intl{Q_\frac 34}|Du|^2\leq C\int_{Q_1}(|u|^3+|P|^\frac 32).$$
\end{lemma}
{\bf Proof}. The estimate is a consequence  of 
(\ref {energy}) by taking proper $\varphi$ and the H\"older inequality.

\begin{lemma}\label{NScompact}
(a) If $\{u, P\}\in L^3(Q_1)\times L^\frac 32(Q_1)$ is a weak solution then $u\in L^{\frac {10}{3}}(Q_\frac 12)$.

(b) If a set $X\subset  L^3(Q_1)\times L^\frac 32(Q_1)$ consists of suitable weak solutions of $\{u, P\}$  with  bounded norms in $L^3(Q_1)\times L^\frac 32(Q_1)$, then $\{u: (u, P)\in X\}$  is compact in $L^3(Q_\frac 12)$.
\end{lemma}

{\bf Proof}
From the previous lemma, we see that $u\in L^\infty((-\frac 12, 0]; L^2(B_\frac 12))\cap L^2((-\frac 12, 0), L^6(B_\frac 12))$ from the Sobolev embedding.

Now we interpolate in space and time to obtain $u\in L^p(Q_\frac 12)$ with bounds since the only solution of the index equation:
$$
\frac 1p=\frac \theta{\infty}+\frac{1-\theta}{2}=\frac \theta 2+ \frac{1-\theta}{6}
$$
is $p=\frac {10}3.$

The compactness in $L^\frac 32(Q_\frac 12)$ follows from Aubin-Lions Lemma since $u_t\in L^\frac 32((-\frac 14, 0], W^{-1, 3}(B_\frac 12))$ and then interpolate with (a) to get the compactness in $L^3(Q_\frac12)$.

\section{  Regularity in the linear scale}
In this section, we will prove a regularity theorem using linearization of Navier-Stoke equations.  Accordingly, we will use the scaling of linear equation in the main iterations. 

Suppose $\{u, P\}$ is a suitable weak solutions of (\ref{NS})  in $Q_1$.
Then its linearization is
\begin{equation}
\left\{
\begin{array}{rll}
u_t+V\cdot \nabla u-\Delta u&=\nabla  P , &\mbox{ in } Q_1,\\
\div u&=0,&\mbox{ in } Q_1,
\end{array}\right.
\end{equation}
where $V$ is a constant vector and shall be considered as an approximation of $u$ in $Q_1$.

If  $\{u, P\}$ is a solution in scale $r$ as 
\begin{equation}
\left\{
\begin{array}{rll}
u_t+V\cdot \nabla u-\Delta u+\nabla  P&=0, &\mbox{ in } Q_r,\\
\div u&=0,&\mbox{ in } Q_r,
\end{array}\right.
\end{equation}
and its {\it linear} scaled solution for continuity or $L^\infty$ scale is defined as 
\begin{eqnarray*}u_r(x, t)&=&u(rx, r^2 t),\\
 P_r(rx, r^2 t)&=&rP(rx, r^2 t),
\end{eqnarray*}
 and these scaled solution would solve the following equation in the unit scale:
\begin{equation}
\left\{
\begin{array}{rcl}\label{linearized}
u_t+rV\cdot \nabla u-\Delta u+\nabla  P&= 0, &\mbox{ in } Q_1\\
\div u&=0,&\mbox{ in } Q_1,
\end{array}\right.
\end{equation}
where we have suppressed the subindex for clarity.

Accordingly if we use  the control of $u$ and $P$  in the unit scale by
$$
\sqrt[3]{\intl{Q_1}|u|^3}+\sqrt[\frac32]{\intl{Q_1}|P|^\frac 32}.
$$
Then scaled control in $Q_r$ would be, the cubical functional as
$$
\sqrt[3]{\intl{Q_1}|u_r|^3}+\sqrt[\frac32]{\intl{Q_1}|P_r|^\frac 32}=\sqrt[3]{\intl{Q_r}|u|^3}+r\sqrt[\frac32]{\intl{Q_r}|P|^\frac 32}
$$
Thus we shall define, 
\begin{equation}\label{Q}
C_r[u, P]=\sqrt[3]{\intl{Q_r}|u|^3}+r\sqrt[\frac32]{\intl{Q_r}|P|^\frac 32}.
\end{equation}
The main analytical result of this section is the following, which is in an iterative form for the convenience of iterations.

\begin{proposition}[Key  Iteration in linear scale]\label{key-linear}
There are universal constants $0<\lambda, \epsilon_0<1$  and $\Lambda$ such that if $\{u, P\}$ is a solution of Navier-Stokes equations  with  approximation (\ref {app} ) \begin{equation}\label{app}
\left\{
\begin{array}{rcl}
u_t+(rV+r u)\cdot \nabla u-\Delta u+\nabla  P&=& 0, \mbox{ in } Q_1\\
\div u&=&0,\mbox{ in } Q_1,
\end{array}\right.
\end{equation}
 with  $|V|\leq 1, 0\leq r\leq 1$ and  $C_1[u, P]\leq \epsilon_0$,

then there is a vector $V_1$ in $R^3$ and a function of time $f(t)$ such that   $|V_1|\leq \Lambda,   C_1[u, P]\leq \frac 12$ and 
\begin{equation}\label{keyinlinear}
C_\lambda[u-V_1, P-f(t)]\leq \frac 12 C_1[u, P].
\end{equation}
\end{proposition}
An easy iteration of the above proposition yields the following.
\begin{theorem}\label{Q}
Let the constants $\lambda, \epsilon_0$ and $ \Lambda$ be as set forth in Proposition \ref{key-linear} and suppose that  $\{u, P\}$ is a solution of (\ref{NS})    in $Q_1$ 

Then  (a)
the smallness condition
\begin{equation}\label{small}
C_1[u, P]\leq {\epsilon_0},
\end{equation} 
 implies that for all $k=1, 2, \dots$, there are vectors $V_k$ and function $f_k(t)$ of time only  so that
$$
C_{\lambda^k}[u-V_k, P-f_k(t)]\leq \frac 1{2^k} C_1[u, P],
$$
and that $V_0=0$ and that $|V_{k}-V_{k-1}|\leq \frac \Lambda {2^{k-1}}C_1[u,P]\leq \frac 1{2^k}$.

(b). $u$ is $C^\alpha(Q_\frac 12)$  provided
\begin{equation}\label{small}
C_1[u, P]\leq \frac{\epsilon_0}{2^5}.
\end{equation} 
\end{theorem}
{\bf Proof.} The proof of (a) is an iteration of (\ref{keyinlinear}). and (b) is a consequence  of the  parabolic Campanato embedding theorem  from the estimates of (a) applied to $u(\frac {x-y}2, \frac{t-s}4)$ for $(y, s)\in Q_\frac 12$ .

We  prove the key linear iteration in the next two  lemmas.

\begin{lemma}[Compact without Compactness in Pressure]\label{app-linear}
For any $\epsilon>0$ there is a $\delta>0$ so that if $\{u, P\}$ is a weak solution of (\ref{app})  for some $|V|\leq M_0$ and $r\leq 1$ with $C_1[u,P]\leq \delta$, then
there are solutions   $v, q$ of  the linearized NS system and $H(x, t)$ of harmonic flow as
\begin{equation}
\left\{
\begin{array}{rcl}\label{linearizedh}
v_t+rV\cdot \nabla v-\Delta v+\nabla  q&=& 0, \mbox{ in } Q_\frac 12\\
\div\, v&=&0,\mbox{ in } Q_\frac12,\\
\Delta H&=&0\mbox{ in } Q_\frac12\\
\end{array}\right.
\end{equation}
such that $C_\frac 12[v, q]\leq 2\delta$, $C_\frac 12[0, H]\leq 2\delta$ and that
\begin{equation}\label{compact}
C_{\frac 12}[u-v, P-H]\leq \epsilon C_1[u, P].
\end{equation}
\end{lemma}
{\bf Proof.} We prove (\ref{compact}) by contradiction.
Suppose it were not true, then there would be an $\epsilon_0>0$ such that for any $n$, there are $u_n, P_n$ with  $C_1[u_n, P_n]\leq \frac 1n$ but 
$$
C_{\frac 12}[u_n-v, P_n-H]\geq \epsilon_0 C_1[u_n, P_n]
$$
for any  solutions $\{v, q, H\}$  of (\ref{linearizedh}).

Now let $\widetilde{u_n}=\frac {u_n}{{C_1[u_n, P_n]}}$ and $\widetilde{P_n}=\frac {P_n}{{C_1[u_n, P_n]}}$.
Then,  we have $C_1[\widetilde{u_n},\widetilde{P_n}]=1$,  $\{ \widetilde{u_n},\widetilde{P_n}\}$ is a solution of
\begin{equation}
\left\{
\begin{array}{rcl}\label{normal}
\widetilde {u_n}_t+(rV+rC_1[u_n, P_n]\widetilde{u_n})\cdot \nabla \widetilde {u_n}-\Delta \widetilde {u_n}+\nabla  \widetilde {P_n} &=0, &\mbox{ in } Q_1\\
\div \widetilde {u_n}&=0,&\mbox{ in } Q_1,
\end{array}\right.
\end{equation}
 and that
\begin{equation}\label{contra}
C_\frac 12[\widetilde u_n-v,P_n-H]\geq \epsilon_0,
\end{equation}
for any solution $\{v, q, H\}$ of (\ref{linearizedh}) with $C_\frac 12[v, q]\leq 2$ and  $C_\frac 12[0, H]\leq 2$.

From Lemma \ref{NScompact},  we have $\widetilde{u_n}$ is bounded $ L^\frac {10}3(Q_\frac 34)$ and  compact in $L^3(Q_\frac 34)$ by the Aubin-Lions Lemma. Therefore we may assume, without lose of generality that $\widetilde u_n\to u_\infty$ in $L^3(Q_\frac 34)$ for some $u_\infty\in L^3$ and that $\widetilde{P_n}\weakto P_\infty$ weakly in $L^\frac 32(Q_\frac 34)$ as well as $u_n\weakto u_\infty$  weakly in $L^2(Q_\frac34)\cap L^2(H^1)$.

Taking divergence in (\ref{normal}), we have
$$
\Delta \widetilde{P_n}=-rC_1[u_n, P_n]\partial_{x_i } \widetilde {u}^j_{n}\partial_{x_j} \widetilde{u}_n^i .
$$ 
Thus we consider
$$\left\{\begin{array}{rll}
\Delta R_n&=-rC_1[u_n, P_n]\partial_{x_i } \widetilde {u}_n^j\partial_{x_j} \widetilde{u}^i_n , &\mbox{ in }Q_\frac 34,\\[2mm]
R_n&=0, &\mbox { on } \partial B_\frac 34\times(-\frac 9{16}, 0].\end{array}\right.
$$
It is easy to see that $R_n\to 0$ in $L^\frac 32(Q_\frac 34)$ strongly and in particular, $C_\frac12 [0, R_n]\leq 1$ for $n$ large.

Now let $H_n=\widetilde{P_n}-R_n$.  Then  $H_n$ is harmonic in  $Q_\frac 34$ and with estimates $C_\frac 12[0, H_n]\leq 2$ for $n$ large.
Now we take weak limit  in (\ref{normal})  as $n\to \infty$ to obtain that
$$
\left\{
\begin{array}{rcl}\label{linearized2}
 (u_\infty)_t+rV_\infty\cdot \nabla  {u_\infty}-\Delta  {u_\infty}+\nabla   {P_\infty}&=& 0, \mbox{ in } Q_\frac34
\\[2mm]
\div  {u_\infty}&=&0,\mbox{ in } Q_\frac 34,
\end{array}\right.
$$
which is a contradiction with (\ref{contra}) if we take $v={u_\infty}$and $H=H_n$ for $n$ large.

\parskip=6pt

{\bf Remark:} We remark that the pressure $P$ in the Navier-Stokes equation lacks compactness since  $P+f(t)$ satisfies the same equation for arbitrary functions $f(t)$. Thus  any compactness statement shall take this into account.  The harmonic flow $H$ in the  above compactness lemma  is a reflection of this fact. This is the reason that we call this lemma {\it Compact without compactness in the pressure}. As matter of fact,  the sequence $f_n(t)$ constructed in  Theorem {\ref{Q}  are not convergent nor under controlled  in  $L^\frac 32$. The idea shall work for problems in other situations.

\parskip=6pt

\noindent{\bf Proof of Proposition \ref{key-linear}.} Let $C_1[u,P]=\delta$. 
From the conditions for ${v, q, H}$, we have  $C_{\frac12}[v, q]\leq 2\delta$, $C_{\frac 12}[0, H]\leq 2\delta$. From the equation for $v$ and $q$,  $v$ is  locally smooth in space and  $C^\frac 13$ in the time variable. Further with the gradient estimates for harmonic function $H$  on each time slice,  and  for each  $(x, t)\in Q_r$, 
\begin{eqnarray}
|v(x, t)-v(0,0)|^3&\leq& C\delta^3 r,\label{v}\\
|H(x, t)-H(0, t)|^\frac 32&\leq & C r^\frac 32\|H(\cdot, t)\|_{L^\frac 32(B_\frac 12)}^\frac 32,\label{H}\\
|v(0, 0)|&\leq & C\delta.
\end{eqnarray}
Now we let $V_1=v(0,0)$ and $f(t)=H(0, t)$ and  integrate (\ref{v}, \ref{H}) in $Q_r$:
\begin{eqnarray*}
\intl{Q_r}|v-V_1|^3&\leq &Cr\delta^3
\\
r^\frac 32\intl{Q_r}|H(x, t)-f(t)|^\frac 32&\leq &
\frac {Cr^3}{r^5}\int_{-r^2}^0r^3\int_{B_\frac 12}|H|^\frac 32dx dt\leq  {C}{r}\intl{Q_\frac 12}|H|^\frac 32\leq Cr\delta^\frac 32.
\end{eqnarray*}
Therefore, 
\begin{eqnarray}
\sqrt[3]{\intl{Q_r}|u-V_1|^3}&\leq&\sqrt[3]{\intl{Q_r}|u-h|^3}+\sqrt[3]{\intl{Q_r}|v-v(0,0)|^3}\notag\\
&\leq&\frac {\epsilon}{r^\frac 53}\sqrt[3]{\intl{Q_\frac 12}|u-h|^3}+\sqrt[3]{\intl{Q_r}|v-v(0,0)|^3}\notag\\
&\leq &\frac {\epsilon}{r^\frac 53}C_1[u, P]+Cr^\frac 13C_1[u,p],\label{u-v}
\end{eqnarray}
and that
\begin{eqnarray}
r \sqrt[\frac 32]{\intl{Q_r}|P-f(t)|^\frac 32}&\leq&r\sqrt[\frac 32]{\intl{Q_\lambda}|P-H|^\frac 32}+r\sqrt[\frac32]{\intl{Q_r}|H-H(0,t)|^\frac 32}
\notag\\
&\leq& \frac {\epsilon}{r^\frac 53}C_1[u, P]+C{r} C_1[u,P].\label{P-H}
\end{eqnarray}
Now we take sum of (\ref{u-v}, \ref{P-H})  to have
$$
C_\lambda[u-V_1, P-f(t)]\leq \frac 12 C_1[u,P], 
$$
if we  first take $r=
\lambda<1$ small enough and then $\epsilon=\epsilon_0.$

This completes the linear approximation theory and the proof of Theorem \ref{Q}.

\section{Nonlinear regularity with small energy}
The goal of this section is to prove that if the following inequality holds in all small scales 
\begin{equation}E_r=E_r(u):=r^4\intl{Q_r}|Du|^2 dxdt=\frac 1{r|B_1|}\int_{Q_r}|Du|^2 dxdt\leq \epsilon,
\end{equation}
then the smallness condition (\ref{small}) holds near $0$ at some small scale,  and thus $u$ is regular near $0$. We have to use nonlinear scaling for Navier Stokes equation which we will review now.

If $\{u, p\}$ is a solution of Navier Stokes equation in $Q_r$, then $u_r(x, t)=r u(rx, r^2 t), P_r(x, t)=r^2P(rx, r^2 t)$ solve the Navier Stokes equation in $Q_1$.

If we use the control on $u, P$ as before, that is
$$
C_1[u,P]=\sqrt[3]{\intl{Q_1}|u|^3}+\sqrt[\frac32]{\intl{Q_1}|P|^\frac 32},
$$
then its {\it nonlinear} scaled control  in $Q_r$ is 
%
then
$$
C_1[u_r, P_r]=r\left(\sqrt[3]{\intl{Q_r}|u|^3}+r\sqrt[\frac32]{\intl{Q_r}|P|^\frac 32}\right).
$$ 
Thus we introduce
\begin{equation}\label{N}
N_r:=N_r[u, P]=rC_r[u, P]=r\left(\sqrt[3]{\intl{Q_r}|u|^3}+r\sqrt[\frac32]{\intl{Q_r}|P|^\frac 32}\right).
\end{equation}
 We also remark that 
 $$
 E_r=E_1(u_r),
 $$
 i.e. $E_r$ is the correct local energy under the nonlinear scale.

\begin{lemma}\label{u-p}
Let $\{u,P\}$ be a suitable weak solution  of (\ref{NS}) in $Q_1$.
Then for $r\leq \frac 12$,
\begin{eqnarray}
(a).  &&\displaystyle r^3\intl{Q_r}|u|^3\leq C(r+\frac {E_1}{r^2}) N_1[u, P]^3+C(\frac 1{r^5}E_1^2+\frac 1{r^3}E_1^\frac 43),\label{A}\\
(b).  && \displaystyle {r^3}\intl{Q_r}|P|^\frac32\leq  C(r+\frac {E_1}{r^2}) N_1[u,P]^{\frac 32}+\frac C{r^5}E_1^2.
\label{B}\end{eqnarray}
\end{lemma}

{\bf Proof.}
We prove (\ref{B}) first. 
Let $\displaystyle\bar{u^i}=\intl{B_\frac34}u^i$. Consider
$$
\Delta P=u^i_{x_j}u^j_{x_i}=((u^i-\bar{u^i})(u^j-\bar u^j))_{x_ix_j}
$$
then apply (b) in Lemma \ref{inmono} to $P$, with $q=3, p=\frac 32$,
to obtain that
\begin{eqnarray*}
\intl{B_r}|P|^\frac 32&\leq &C\intl{B_1}|P|^\frac 32+\frac {C}{r^3}\int_{B_\frac 34}|u-\bar u|^3\\
&\leq&C\intl{B_1}|P|^\frac 32+\frac {C}{r^3}\left(\int_{B_\frac 34}|u-\bar u|^2\right)^\frac 34\left(\intl{B_\frac34}|u-\bar u|^6\right)^\frac 14\\
&\leq&C\intl{B_1}|P|^\frac 32+\frac {C}{r^3}\left(\int_{B_\frac 34}|u-\bar u|^2\right)^\frac 34\left(\intl{B_1}|Du|^2\right)^\frac 14\\
&\leq&C\intl{B_1}|P|^\frac 32+\frac {C}{r^3}\left(\int_{B_\frac 34}|u-\bar u|^2\right)^\frac 12\intl{B_1}|Du|^2\\
&\leq&C\intl{B_1}|P|^\frac 32+\frac {C}{r^3}\sup_{t\in (-\frac 14, 0]} \left(\int_{B_\frac 34}|u|^2\right)^\frac 12\intl{B_1}|Du|^2.
\end{eqnarray*}

Integrate  in $t$ over $(-r^2, 0]\subset(-\frac 14, 0]$: 
\begin{eqnarray*}
\frac 1{r^3|B_1|}\int_{Q_r}|P|^\frac 32&\leq &C\intl{Q_1}|P|^\frac 32+\frac {C}{r^3}\sup_{t\in (-\frac 14, 0]} \left(\int_{B_\frac 34}|u|^2\right)^\frac 12\int_{Q_1}|Du|^2.\\
\end{eqnarray*}

Now we invoke the local energy estimates in  Lemma \ref{local-energy},  and  the elementary inequalities  $\sqrt{A+B}\leq \sqrt A +\sqrt B$ and $2AB\leq A^2+B^2 $ to obtain,
\begin{eqnarray*}
\frac 1{r^3|B_1|}\int_{Q_r}|P|^\frac 32&\leq &C\intl{Q_1}|P|^\frac 32+\frac {C}{r^3}\sqrt{\int_{Q_1}|u|^3+\int_{Q_1}|P|^\frac 32}\intl{Q_1}|Du|^2\\
&\leq &CN_1[u,P]^\frac 32+\frac {C}{r^3}\left(\sqrt{N_1[u,P]^3+N_1[u,P]^\frac 32}\right)E_1\\
&\leq& CN_1[u,P]^\frac 32+\frac {C}{r^3}\left(N_1[u,P]^\frac 32E_1+N_1[u,P]^{\frac 34}E_1\right)\\
&\leq&CN_1[u,P]^\frac 32+\frac {C}{r^3}N_1[u,P]^\frac 32E_1+N_1[u,P]^{\frac 32}+\frac {C^2}{r^6}E_1^2
\end{eqnarray*}
 which is (\ref{B}) after multiplying by $r$.

The proof for (\ref{A}) is similar and  is actually easier. Recall  (\ref{ininter}):
$$
 \intl{B_r}|u|^3 dx\leq C\intl{B_1}|u|^3 dx+\frac {C}{r^3}\left(\int_{B_\frac 34}|u|^2\right)^\frac 12\int_{B_1} |Du|^2.
$$
Integrate  in $t$ over $(-r^2, 0]\subset(-\frac 14, 0]$: 
\begin{eqnarray*}
\frac 1{r^3|B_1|}\int_{Q_r}|u|^3&\leq &C\intl{Q_1}|u|^3+\frac {C}{r^3}\sup_{t\in (-\frac 14, 0]} \left(\int_{B_\frac 34}|u|^2\right)^\frac 12\int_{Q_1}|Du|^2\\
&\leq &C\intl{Q_1} |u|^3+\frac {C}{r^3}\sup_{t\in (-\frac 14, 0]} \left(\int_{B_\frac 34}|u|^2\right)^\frac 12\int_{Q_1}|Du|^2.
\end{eqnarray*}

Now  as before with the local energy estimates in  Lemma \ref{local-energy},  and the elementary inequalities $\sqrt{A+B}\leq \sqrt A +\sqrt B$ and H\"order inequalty: $A^\alpha B^{1-\alpha}\leq \alpha A^\frac 1\alpha +(1-\alpha)B^\frac 1{1-\alpha} $ to obtain,
\begin{eqnarray*}
\frac 1{r^3|B_1|}\int_{Q_r}|u|^3&\leq &C\intl{Q_1}|u|^3+\frac {C}{r^3}\sqrt{\int_{Q_1}|u|^3+\int_{Q_1}|P|^\frac 32}\intl{Q_1}|Du|^2\\
&\leq &CN_1[u,P]^3+\frac C{r^3}\left(N_1[u,P]^\frac 32E_1+N_1[u,P]^{\frac 34}E_1\right)\\
&\leq& CN_1[u,P]^3+C\left(\frac 34C_1[u,P]^3+(\frac  1{2r^6}E_1+\frac{3}{4r^4}E_1^{\frac 43})\right),
\end{eqnarray*}
which is  (\ref{A}) after multiplying by $r$.

\begin{lemma}[Key Nonlinear Iteration] \label{epsilonlambda}There are constants $0<\lambda, \epsilon_1<1$
such that if $E_1\leq \epsilon_1$, then 
\begin{equation}\label{nonlinearkey}N_\lambda\leq \frac 12N_1+E_1^\frac 13.\end{equation}
\end{lemma}

{\bf Proof.} Apply  $(A+B+C)^\alpha\leq A^\alpha+B^\alpha+C^\alpha$
to (\ref{A}) for $\alpha=\frac 13$ and to (\ref{B}) for $\alpha=\frac 23$  and then add them to obtain:
\begin{equation}
N_r[u, p]\leq \left(Cr^\frac 13+\frac {CE_1^\frac 13}{r^\frac 23}+Cr^\frac 32+\frac {CE_1^\frac 23}{r^\frac 43}\right)N_1+C(\frac {E_1^\frac 23}{r^\frac {5}3}+\frac {E_1^\frac 49}{r}+\frac {E_1^\frac 43}{r^\frac {10}2}).
\end{equation}
Therefore  we take $r=\lambda $ small so that ${C}\lambda^{\frac 13}+C{\lambda^\frac 23}\leq \frac 14$  and then take $\epsilon$ small so that $C\frac{\epsilon^{\frac 13}}{\lambda^\frac 23}+C\frac{\epsilon^{\frac 23}}{\lambda^\frac 43}\leq \frac 14$ as well as
$C(\frac {\epsilon^\frac 23}{\lambda^\frac {5}9}+\frac {\epsilon^\frac 13}{\lambda}+\frac {\epsilon^{\frac {11}9}}{\lambda^{\frac {10}2}})\leq 1,$
to arrive at (\ref{nonlinearkey}). The lemma is proved.

\begin{lemma}
Let $\epsilon_1, \lambda$ as set forth in Lemma \ref{epsilonlambda} and  if $E_r\leq  \epsilon<\epsilon_1<1$  for all $0<r<1$, then,

(a) for any $r\leq 1$, $N_{r\lambda }\leq \frac 12 N_{r}+\epsilon^\frac 13,$

(b) for any   $k\in \mathbf N$, $N_{\lambda^k}\leq \frac 1{2^k}N_1+2\epsilon^\frac 13.$

\end{lemma}

{\bf Proof.}
(a) is the scaling of the conclusion of the previous lemma and (b) is an iteration of (a).

\begin{proposition} 
If $E_r\leq \epsilon_2:=\min(\frac {\epsilon_0^3}{2^{21}}, \epsilon_1)$, for all $r<1$, then $u$ is H\"older continuous in $Q_{r_0}$ for $r_0=\lambda^{6}(\frac{N_1}{\epsilon_0})^{\frac {ln \lambda}{ln2}}$.
\end{proposition}
{\bf Proof.} Recall $N_1=N_1[u, P]$, and  take $k\geq\log_2(\frac {N_1}{\epsilon_0})+6$ as well as $\epsilon_2\leq \frac{\epsilon_0^3}{2^{21}}$, then
$$
N_{\lambda ^k}\leq \frac {N_1}{2^k}+2\epsilon_2^\frac 13\leq \frac {\epsilon_0}{2^6}+\frac {\epsilon_0}{2^6}=\frac{\epsilon_0}{2^5}.
$$
Now we can apply Theorem \ref{Q} to obtain  $u$ is continuous in $Q_{\frac {\lambda^k}2}$.

\begin{theorem}[Main theorem]
If $E_r\leq \epsilon_2$ for all $0<r\leq R$, then $u$ is regular near the origin. Any suitable weak solution of the Navier Stokes equations not continuous in  a set, then such set has measure 0  with respect to  the one dimensional parabolic Hausdroff measure.
\end{theorem}
{\bf Proof} We notice that $E_\frac rR(u_R)=E_r(u)\leq \epsilon_2$,  for all $\frac rR\leq 1$. Thus $u_R$ is H\"older continuous in $Q_{r_0}$ and so is $u$ in $Q_{\frac{r_0}2}$. The last conclusion from the standard covering argument.

{\bf Acknowledgement} This research is support in part by a grant from Simons Foundation.

\begin{tabular}{ll}
Department of Mathematics&School of Mathematical Sciences \\
University of Iowa&Shanghai Jiaotong University\\
Iowa City, IA52242&Shanghai, 200240\end{tabular}
\end{document}